\documentclass[review]{elsarticle}

\usepackage[linktocpage=true,pagebackref=true,linkcolor=black,citecolor=black]{hyperref}

\journal{Journal of \LaTeX\ Templates}

\usepackage{verbatim}
\usepackage[english]{babel}
\usepackage{amssymb}
\usepackage{amsthm}
\usepackage{amsmath,amssymb,amsfonts,graphicx,float,multicol,fleqn}
\usepackage{subcaption}
\usepackage{multirow}
\usepackage{booktabs}
\usepackage{caption}
\usepackage[flushleft]{threeparttable}
\UseRawInputEncoding

\usepackage{algorithm}
\usepackage[noend]{algpseudocode}

\usepackage[letterpaper, margin=1in]{geometry}
\usepackage[dvipsnames]{xcolor}
\usepackage{enumerate}

\usepackage{caption} 
\theoremstyle{definition}

\begin{document}
\begin{frontmatter}
\title{Designing of knowledge-based potentials via B-spline basis functions for native proteins detection}

\author[1]{Elmira Mirzabeigi
\corref{mycorrespondingauthor}}
\ead{e.mirzabeigi@modares.ac.ir}
\author[2]{Saeed Mortezazadeh}
\ead{s.mortezazadeh@modares.ac.ir}
\author[1]{Rezvan Salehi}
\ead{r.salehi@modares.ac.ir}
\author[2]{Hossein Naderi-Manesh}
\ead{naderman@modares.ac.ir}
\cortext[mycorrespondingauthor]{Corresponding author}
\address[1]{Department of Applied Mathematics,
	Faculty of Mathematical Sciences,
	Tarbiat Modares University,
	Tehran, Iran}
\address[2]{Department of Nanobiotechnology,
	Faculty of Biological Sciences,
	Tarbiat Modares University,
	Tehran, Iran}
\begin{abstract}
Knowledge-based potentials were developed to investigate the differentiation of native structures from their decoy sets. This work presents the construction of two different distance-dependent potential energy functions based on two fundamental assumptions using mathematical modeling. Here, a model was developed using basic mathematical methods, and the carbon-alpha form is the simplest form of protein representation. We aimed to reduce computational volume and distinguish the native structure from the decoy structures. For this purpose, according to Anfinsen’s dogma, we assumed that the energy of each model structure should be more favorable than the corresponding native type. In the second one, we thought that the energy difference between the native and decoy structures changes linearly with the root-mean-square deviation of structures. These knowledge-based potentials are expressed by the B-spline basis functions of the pairwise distances between C$_{\alpha}$-C$_{\alpha}$ of inter-residues. The potential function parameters in the above two approaches were optimized using the linear programming algorithm on a massive collection of Titan-HRD and tested on the remainder. We found that the potential functions produced by Anfinsen’s dogma detect native structures more accurately than those developed by the root-mean-square deviation. Both linear programming knowledge-based potentials (LPKP) successfully detect the native structures from an ensemble of decoys. However, the LPKP of the first approach can correctly identify 130 native structures out of 150 tested cases with an average rank of 1.67. While the second approach LPKP detects 124 native structures from their decoys. We concluded that linear programming is promising for generating knowledge-based potential functions. All the training and testing sets used for this work are available online and downloaded from \url{http://titan.princeton.edu/HRDecoys}.
\end{abstract}
\end{frontmatter}

\section{Introduction}
Proteins are important macro-molecules that are involved in all cellular processes and their functions are directly related to the 3D structure. The prediction of the protein structure is one of the incentive questions in computational biology. Although we know millions of protein sequences, less than a hundred thousand protein structures have been found until now \cite{PlosONE}
. It can be indirectly concluded that determining the protein structure at the atomic level is overwhelming. Reports on \textit{ab initio} technique developments show significant progress in predicting protein structure in previous years, however, the quality of models does not have enough accuracy to be useful for biologists \cite{Zhang}. Hence, the bio-informatics approach is most widely used to predict the tertiary structure of proteins.
\\
According to Anfinsen's dogma, the native structure is determined by the amino acid sequence of protein which means at the environmental condition that folding occurs, the native structure is formed at the global minimum free energy \cite{anfinsen}. Force fields have been developed to calculate the potential energy of molecular systems which refer to the functional form and parameter sets that can be derived from empirical and theoretical studies. The protein structure is considered atomic or coarse-grained so that potential functions can include pairwise interactions, side-chain orientations, secondary structural preferences, solvent exposure, and other geometric properties of proteins \cite{Zhou, Skolnick}. The accuracy of these potential functions can be traced by evaluating the ability to detect the native structure from a decoy ensemble. There are two general knowledge-based and physics-based approaches to introduce potential functions \cite{Summa, Zhu, Chopra, Amautova, Bhattachary}. Knowledge-based potentials are simple potentials extracted from the protein data bank designed to improve the quality of protein models whereas, in physics-based potentials, the chemical properties of the molecule are also taken into account. Various methods have been developed to increase the accuracy of these potential functions in recent years \cite{Rohl, Zhang2003, Benkert, Zhang2004}.
\\
Nowadays, many methods like Rosetta, Alpha-fold, and Michael Feig's refinement work are routinely getting high-accuracy structures. In addition, Baker and Feig's refinement methods are able to generate ensembles for near-native structures and score the best structures in those ensembles to produce improved structures \cite{Akhter, Uziela, Akhter2, Zhang2, ZhangGJ, Jing, Alford, Hou, Wang, Li}. In most cases, knowledge-based potentials are used to sort decoy structures by calculating similarity to the native structure. The inability of the proposed model to identify the native structure as the lowest potential energy is regarded as the error rate of the potential functions. There are two methods to generate knowledge-based potentials. The first method is to use the Boltzmann inversion equation to convert the distribution of geometric properties of known structures to potential energy functions. This method is also called the statistical potential obtained from the ratio of observed frequencies to the reference state. Therefore, several models have been developed for representing reference states such as Sippl's assumption, distance-scaled finite ideal-gas reference state (DFIRE), and discrete optimized protein energy (DOPE). In the second method, the potentials are extracted from a training process to quantitatively differentiate between the incorrectly folded models and native structures \cite{PlosONE}.
\\
Using physical potential functions and considering more than 40 physical constraints in a linear programming problem, Rajgaria et.al \cite{Rajgaria} have achieved better results than previous studies on the Titan-HRD protein data set. Despite the desirable results, it is inevitable to consider several constraints in the optimization problem using the physical potential functions. In some studies, the least-square problem has been used to optimize knowledge-based potentials. Since the least-square problem is based on norm-2, the energy difference between the native structure and the decoy can give both positive and negative values. Therefore, the use of this method in protein modeling makes Anfinsen's dogma not applicable to the potential production process. In these studies, choosing the best decoy instead of recognizing the native structure is the main parameter for model evaluation. The best decoy is determined by calculating the minimum distance and the lowest energy difference with the native structure. In addition to the inability to distinguish the native structure from the decoy, there is also a high computational cost for using the least-squares method \cite{PlosONE, CarlsenProteins}.
\\
We considered the pairwise distances between the C$_{\alpha}$ of amino acids so that the knowledge-based potentials are formulated using B-spline basis functions. We introduce two different optimization processes to obtain potential parameters. Both are based on Anfinson's dogma that in the first method only the energy difference of the native and decoy structure is included in the modeling while in the latter method, it is assumed that the energy difference is linearly proportional to the distance of structures. We solved them by using MATLAB software and checked the results of the two methods. Finally, our results were cross-checked with previously reported data. Our LPKPs have better detection for native structures compared to other methods. Both LPKPs are able to identify 130 and 124 native structures out of 150 tested cases, respectively.
\section{Materials and Methods}
In this model, the amino acids are represented by the location of C$_{\alpha}$ atoms. Hence a protein with $m$ amino acids is expressed by a vector, $C=(c_{1}, c_{2}, ..., c_{m})$ where $c_{k}$ is the location of C$_{\alpha}$ of the $k$-th amino acid having the following three-dimensional coordinates:
\begin{equation}
\hspace*{\fill} c_{k} = (x_{k},y_{k},z_{k}). \hspace*{\fill}
\end{equation}
We denote the native structure of each protein with $N=(n_{1}, n_{2}, ..., n_{m})$ and decoy structures with $D_{i}=(d_{1}, d_{2}, ..., d_{m})$, where $i$ is the decoy number.
\subsection{Geometric Distance Between Two Protein Structures}
The Euclidean distance between two points is the basis for calculating the geometric distance between two protein structures. Root-mean-square deviation (RMSD) method was used for this purpose so that the first protein structure is considered to be the reference and the second structure is best fitted to the first one using a transformation operator $G$:
\begin{equation}
	\hspace*{\fill} RMSD = \min_{G}\sqrt{\frac{\sum_{k=1}^{m}\|n_{k}-G(d_{k})\|^{2}}{m}}, \hspace*{\fill}
\end{equation}
where $n_{k}$ and $d_{k}$ are the $k$-th location of the C$_{\alpha}$ atoms in native reference structure and decoy structures, respectively, and $m$ is the number of amino acids. The transformation $G$ involves a set of translation and rotation operations to obtain the best structural alignment between two protein structures and it does not change the shape or size of the proteins \cite{McLachlan, Horn, Coutsias}. RMSD as defined above is a norm-2 metric distance \cite{Kaindl}. The norm-2 (also written "L$^{2}$-norm") $\|x\|$ is a vector norm defined for a vector $x^{T} = \left[x_{1}, x_{2}, \ldots x_{n}\right]$ by
\begin{equation}
	\hspace*{\fill} \|x\| = \sqrt{\sum^{n}_{k=1}|x_{k}|^{2}}, \hspace*{\fill}
\end{equation}
where $|x_{k}|$ denotes the absolute value of $x_k$ \cite{Johnson}. 
\subsection{Knowledge-based Potential Function}
According to Anfinsen's dogma, potential energy functions must be obtained in such a way that the energy of the protein's native structure has the lowest value in a set of decoy structures \cite{anfinsen}. This hypothesis is shown in the following constraint:
\begin{equation}\label{eq:1}
	\hspace*{\fill} E_{D_{i}} - E_{N} > \varepsilon, \hspace*{\fill}
\end{equation}
in this equation, $E_{N}$ is the energy vector of a native structure and $E_{D_{i}}$ is the energy vector of the $i$-th decoy structure of the relevant protein \cite{Rajgaria}. If we have $p$ proteins and each protein has $i$ decoys, then constraint \ref{eq:1} is reformulated as below:
\begin{equation}\label{eq:2}
	\hspace*{\fill} \sum_{p,i}(E_{p,i}(X) - E_{p,n}(X))> \varepsilon, \hspace*{\fill}
\end{equation}
where $X$ is a vector of parameters that indicates the number of interactions between different types of amino acids at different distances of C$_{\alpha}$-C$_{\alpha}$ in the protein structure. Given that the number of natural amino acids is 20, the total number of types of interactions is equal to 210. Also, each type of interaction is expressed with eight parameters representing potential energy at different intervals. Therefore, the total number of parameters needed to calculate the potential energy between different types of amino acids at different distances is equal to $210 \times 8 = 1680$. In equation \ref{eq:2}, $\varepsilon$ is an arbitrary number in which different values have been tested and the best of them has been 0.01. In the $p$-th protein, $E_{p, i}$ is the energy matrix of the $i$-th decoy, and $E_{p, n}$ is the energy matrix of the native. Now equation \ref{eq:2} can be rewritten as follow, where parameters $S_{p}$ are positive frail variables:
\begin{equation}\label{eq:3}
	\hspace*{\fill} (E_{p,i} - E_{p,n} ) X - S_{p} \geq\varepsilon. \hspace*{\fill}
\end{equation}
In the first approach, the knowledge-based potential functions were optimized through the criteria formulated in equation \ref{eq:3}. Also, assuming that the energy difference between the native and the decoy structure can be related to the geometric distance between them, a similar condition can be extended to optimize the potential functions \cite{PlosONE, CarlsenProteins, Rogen}. This constraint is shown below:
\begin{equation}\label{eq:4}
	\hspace*{\fill} E_{D_{i}} - E_{N} \propto Dis_{D_{i},N}, \hspace*{\fill}
\end{equation}
where, $Dis_{D_{i},N}$ is the distance between decoy $D_{i}$ and native $N$. The equation \ref{eq:4} is rewritten in the following:
\begin{equation}\label{eq:5}
	\hspace*{\fill} E_{D_{i}} - E_{N} = \alpha . Dis_{D_{i},N}. \hspace*{\fill}
\end{equation}
In this equation, $\alpha$ is constant that is considered as 1 in \cite{Rogen}, so the constraint \ref{eq:5} is rewritten for $p$ proteins as below, where $D_{p, i, n}$ is the distance between the native structure and $i$-th decoy structure of $p$-th protein.
\begin{equation}\label{eq:6}
	\hspace*{\fill} \sum_{p,i} E_{p,i}(X) - E_{p,n}(X) \leq D_{p,i,n}. \hspace*{\fill}
\end{equation}
The extended uniform cubic B-spline possesses the convex hull property, symmetry, and geometric invariability \cite{Hollig}. These features are convincing to use this function as the basis for the potential functions. We used this function at eight uniform intervals for each of the 210 types of interactions between amino acids shown in Table 1. So, the energy of all decoy and native structures is calculated using the equation below \cite{PlosONE, Rogen}:
\begin{equation}\label{eq:7}
	\hspace*{\fill} E(\theta) = \sum_{i<j}\sum_{p} X^{aa(i),aa(j)}_{p}B_{p}(r_{i,j}). \hspace*{\fill}
\end{equation}
In this equation, $aa(i)\in \{1, . . . ,20\}$ is the amino acid type of the $i$-th C$_{\alpha}$ and $B_{p}(r_{i,j})$ is the $p$-th B-spline basis function evaluated on the distance between the $i$-th and $j$-th C$_{\alpha}$, shown in Figure 1. Also, $X^{aa(i), aa(j)}_{p}$ are the model parameters determined by the optimization according to the following.

\begin{minipage}{\textwidth}
	\begin{minipage}[h]{0.47\textwidth}
		\centering
		\includegraphics[width=\textwidth]{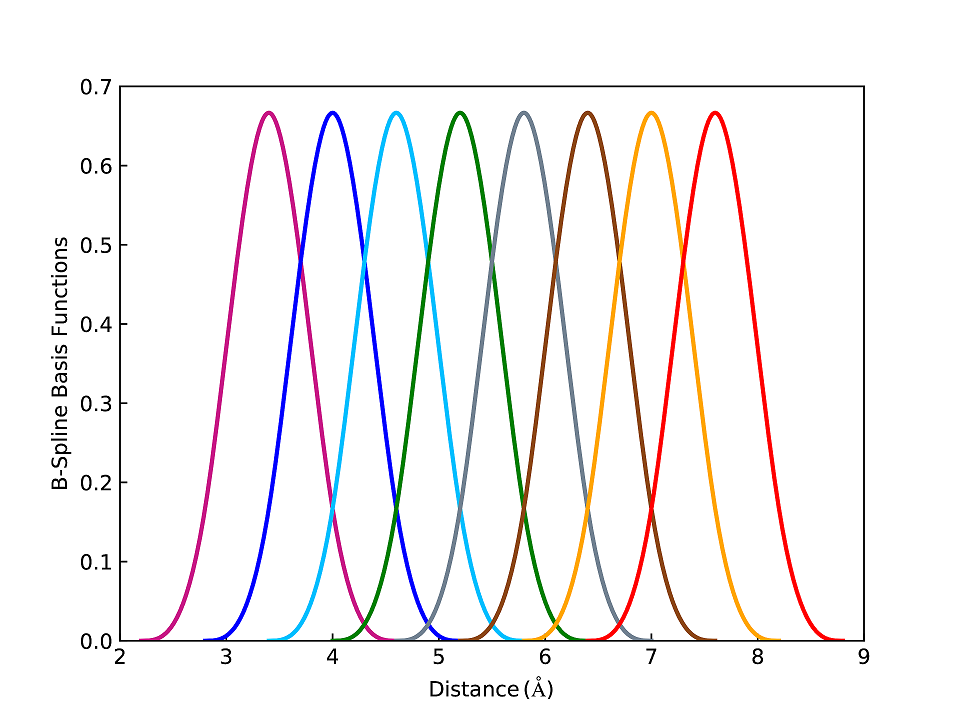}\label{fig1}
		\captionof{figure}{The eight cubic B-spline basis functions $B_{1}$,$...$, $B_{8}$ used in the knowledge-based potentials.}
	\end{minipage}
	\hspace{0.1\textwidth}
	\begin{minipage}[h]{0.4\textwidth}
		\centering
		\captionof{table}{The bin width distances for eight B-spline basis functions.}
		\resizebox{.4 \textwidth}{!}{
		\begin{tabular}{cc}\hline
			\textbf{ID} & \textbf{$C_{\alpha}$ Distance (\AA)}\\\hline\hline 
			1 & 2.2 - 4.6\\
			2 & 2.8 - 5.2\\
			3 & 3.4 - 5.8\\
			4 & 4.0 - 6.4\\
			5 & 4.6 - 7.0\\
			6 & 5.2 - 7.6\\
			7 & 5.8 - 8.2\\
			8 & 6.4 - 8.8 \\
		\end{tabular}
	}
	\end{minipage}
\end{minipage}

\subsection{Optimization Problems}
We designed the energy function using an optimization procedure, assuming that for every native structure $N_{n,p}$, we have a set of decoy structures $D_{i,p}$. At each stage, $(E_{i,p} - E_{n,p})$ is the energy difference, and $D_{i,n,p}$ is the corresponding distance between the native structure and $j$-th decoy in the $p$-th protein. Also, the sum of the frail variables constraint Eq.\ref{eq:2} was minimized in each scheme. Since the frail variables are positive, the condition $S_{p} \geq 0$ was added to each optimization problem. For vector $X$, we consider the subscription of general constraints \cite{Rajgaria} for the condition $-4 \leq X \leq 4$.
\subsubsection{First optimization scheme}
The first optimization problem can be written into linear programming which used Anfinsen's dogma as below, which was named with LPKP$^{1}$:
\begin{align}\label{eq:9}
		\min_{X , S_{1}, ..., S_{p}} \qquad & \sum S_{p} , \\
		s.t \qquad& (E_{i,p} - E_{n,p} ) X - S_{p} \geq\varepsilon ,\\
		& S_{p} \geq 0, \\
		& -4 \leq X \leq 4 .\\
\end{align}
where the vector $X$ is a set of parameters of the energy functions, and $p$ is the number of proteins in the training set.
\subsubsection{Second optimization scheme}
The second optimization problem (called LPKP$^{2}$) can be written into linear programming using Eq.\ref{eq:9} and the relationship between the distance and energy difference of the two structures is as follows:
\begin{align}\label{eq:8}
		\min_{X , S_{1}, ..., S_{p}} \qquad & \sum S_{p} , \\
		s.t \qquad& (E_{i,p} - E_{n,p} ) X - S_{p} \geq\varepsilon ,\\
		& (E_{i,p} - E_{n,p} ) X \leq D_{i,n,p} ,\\
		& S_{p} \geq 0, \\
		& -4 \leq X \leq 4 .\\
\end{align}
This scheme is similar to the first scheme in terms of protein number, optimization method, and the number of parameters $X$ except for the constraint \ref{eq:5}.
\subsection{Training and Test Sets}
The high-resolution decoy set contained 1400 protein structures, with 500-1600 decoys for each protein \cite{Rajgaria}, but we deleted proteins in which decoys and native structures do not have the same number of amino acids. Finally, the 1370 proteins remained from the Titan-HRD. The entire set of protein decoy structures has been made available at \url{http://titan.princeton.edu/HRDecoys/}. From 1370 proteins used for decoy generation, 1220 proteins were randomly selected for training processes. Since each protein has at least 500 decoy structures, therefore, all decoys of each protein were sorted based on their C$_{\alpha}$ RMSD, and then 500 decoys were randomly selected to cover the whole RMSD range. This arrangement of a training set has 500$\times$ 1220 = 610000 decoy structures. Because of computational limitations, it is not possible to include all 610000 decoys in the training step. Therefore, we reduced 500 to 45 decoys per protein to set up the training procedure with 60000 decoy structures. The remaining 150 proteins of the Titan-HRD set were used for testing the obtained potential functions. In this step, also 500 decoys were selected using the same technique explained above in generating the training set. A summary of the test and training data sets is given in Table \ref{table2}.
\begin{table}[h]
	\centering
	\caption{Properties of the different protein decoy sets used in this study.}\label{table2}
	\begin{threeparttable}
		\begin{tabular}{p{3cm}cccc}
			\textbf{Decoy set} & \textbf{Nprot $^{a}$} & \textbf{Nres $^{b}$} & \textbf{Ndecoys $^{c}$} & \textbf{RMSD}\\\hline\hline
			Titan-HRD & 1220 & 111.80 & 500 & 2.46\\
			Titan-HRD$^{*}$ & 150 & 103.95 & 500 & 2.51\\
			\hline 
		\end{tabular}
		\begin{tablenotes}
			\item Training set (Titan-HRD) and test set (Titan-HRD$^{*}$) from the Titan high resolution decoy set at http://titan.princeton.edu/HRDecoys/.
			\item RMSD is the distance measured between the decoys and the corresponding native structures averaged over all decoys and all proteins.
			\item $^{a}$ Number of proteins in each set.
			\item $^{b}$ The average number of residues in each set.
			\item $^{c}$ The average number of decoys in each set.
		\end{tablenotes}
	\end{threeparttable}
\end{table}

\section{Results and Discussion}
The existence of effective factors has led to the use of various methods for the production of knowledge-based potential functions. In this study, we have modeled these functions taking into account a limited number of these parameters. The knowledge-based potential functions are designed according to eight uniform cubic B-spline functions and the distances between C$_{\alpha}$ atoms in protein structures. These potential functions were extracted from the information of 60000 decoy structures. Two sets of energy parameters were constructed based on two different optimization schemes of LPKP$^{1}$ and LPKP$^{2}$. The main constraint in the LPKP$^{1}$ model is Anfinsen's dogma \cite{anfinsen}. Hence, the energy difference between the decoy and the native structure must always be greater than a positive constant value called $\varepsilon$, which is shown in constraint \ref{eq:3}. In the LPKP$^{2}$ model, the effect of RMSD is examined on our scheme by adding constraint \ref{eq:6}. The objective function of each scheme has minimized the sum of the frail variables, discussed in the methods section. The zero value for an objective function means that there are no infringements. In this study, the objective function has obtained a value close to zero, which represents the very low error of our introduced optimization model. Previously, the optimization schemes \ref{eq:5} and \ref{eq:8} were described, and the method was generally explained in Algorithm \ref{algorithm1}.

\begin{algorithm}[bt]
	\caption{Execution of knowledge-based potential function}\label{algorithm1}
	\textbf{Input:} M, $N_{decoys}$.\\
	\Comment{$M:$ Number of Proteins}.\\
	\Comment{$N_{decoys}:$ Number of Decoys in a Similar Protein}.\\
	\textbf{Output:} $E_{i,j}$, $D_{i,n,p}$.\\
	\Comment{$E_{i,j}:$ Energy difference between two structures.}\\
	\Comment{$D_{i,n,p}:$ RMSD distance between two structures.}\\
	\Comment{$\theta^{1}:$ The model parameters determined by the optimization \ref{eq:5}}.\\
	\Comment{$\theta^{2}:$ The model parameters determined by the optimization \ref{eq:8}}.\\
	\begin{algorithmic}[1]
		\For {$i=1$ to $M$}
		\State Read $C_{\alpha}$ three-dimensional coordinates for native structure.
		\State Solve $Dis_{i,1}$ for each pair of $C_{\alpha}$ in native structure.
		\State Solve problem $E_{Native_{i}}$ \ref{eq:7} for native structure.
		\For {$j=1$ to $N_{decoys}$}
		\State Read $C_{\alpha}$ three-dimensional coordinates for decoy structure.
		\State Solve $Dis_{i,j+1}$ for each pair of $C_{\alpha}$ in decoy structure.
		\State Solve problem $E_{Decoy_{i,j}}$ \ref{eq:7} for decoy structure.
		\State Set $E_{i,j} = E_{Decoy_{i,j}} - E_{Native_{i}}$.
		\State Solve RMSD problem $D_{i,n,p}$ \ref{eq:9}.
		\EndFor
		\EndFor
		\State Solve problem \ref{eq:5} for $\theta^{1}$.
		\State Solve problem \ref{eq:8} for $\theta^{2}$.
	\end{algorithmic}
\end{algorithm}

At first, the value of $\varepsilon$ considered 0.01 in schemes \ref{eq:5} and \ref{eq:8}. The optimization schemes were implemented with $\varepsilon$ and approximately 60000 structures. Since the ability to identify between the native and native-like structures is a considerable standard for any potential function, the optimization results have been tested to detect native structures. The LPKP examined 500 decoys for each of the 150 test proteins of the Titan-HRD decoy set. In this examination, the comparative rank of the native among decoy structures has been calculated. An ideal potential function should be able to detect rank 1 for the native structures of all the proteins in the test set. As we know, it is an accepted condition that the test set should not share with the training set since it invalidates the potential energy appeasement. Hence, our test set was carefully selected so that the training and test sets had no overlap. The examination results are presented in Table \ref{table4}. According to that, 130 native structures were correctly recognized among 150 proteins with an average rank of 1.67 in method \ref{eq:5}. Moreover, native structures were distinguished in 124 proteins of scheme \ref{eq:8} with an average rank of 2.83. The results represent that our proposed method is significantly more accurate than previous methods \cite{Rajgaria, LKF, TE13, HL}.
\begin{table}[h]
	\centering
	\resizebox{.8 \textwidth}{!}{
		\begin{threeparttable}
			\caption{Rankings of the native structures using the LPKP$^{1}$ and LPKP$^{2}$ on Titan-HRD test set.}
			\label{table4}
			\begin{tabular}{cccccccccccc}
				\textbf{ID} & \textbf{LPKP$^{1}$} & \textbf{LPKP$^{2}$} &  \textbf{ID} & \textbf{LPKP$^{1}$} & \textbf{LPKP$^{2}$} & \textbf{ID} & \textbf{LPKP$^{1}$} & \textbf{LPKP$^{2}$} &  \textbf{ID} & \textbf{LPKP$^{1}$} & \textbf{LPKP$^{2}$}\\\hline\hline
				1em9A &1& 1 & 1fb1A & 1 & 1 & 1chd- & 1 & 1 & 1faq- & 1 & 1 \\
				1iioA  &1& 1 & 1hf9A & 1 & 2 & 2drpA & 1 & 1 & 1b0yA & 1 & 1\\
				1qqvA &1& 1 & 1b1bA & 1 & 1 & 1ci5A & 1 & 1 & 1tmy- & 1 & 2\\
				1b9lA	&1&	1 & 1bik- & 1 & 1 & 1k3bB & 1 &	1 & 1gd5A	& 1 & 1\\
				1ap0- &1& 1 & 1hgvA	& 39 & 35 & 1fxkC & 1 & 1 & 1hd0A & 1 &	1\\
				1ghc- &1& 1 & 1tyfA & 1 &	1 & 2ech- & 1 & 1 & 1aplD	& 1 & 1\\
				1u2fA	&1&	1&1ag4-	&1 &1 & 1af8-	&1&	1&3monB	&1 &4\\
				1eptB &2&	1&1fe6A	&1&	1&1cjgA	&1 &12 & 1quqB	&1&	2\\
				1aq3A	&1&	4&1b2iA	&2&	1&1g10A&4 &1 & 1eqiA	&1&	1\\
				1c3kA &1&	1&1cm0B	&1 &1 & 1n72A	&1&	1&1cl3A	&1 &1\\
				1l0oC	&1&	2&1dpuA	&1&	1&1k5yR &1 &23 & 1imt-	&1&	1\\
				1etpA	&1&	1&1b2pA	&1&	1&1eal-	&1 &1 & 1rof-	&1&	1\\
				1jq4A	&1&	1&1ahjA	&1 &1 & 1jpyA	&1&	1&1cmaA	&1 &6\\
				1a10I	&1&	1&1hjrA	&1&	1&1qc7A	&1 &1 & 1d7bA	&1&	1\\
				1lcl-	&1&	1&1hks-	&1&	1&1bccH	&1 &3 & 1a14H	&1&	1\\
				1hlb-	&1&	1&1ec5A	&1 &1 & 1ndoB	&1&	1&1qckA	&1 &1\\
				1gnf-	&1&	1&1c6vX	&1&	1&1k99A	&9 &1 &1irqA	&3&	2\\
				1lghB	&1&	39&1c7uA	&3&	1&1i7kA	&1 &1 & 1scjB	&1&	1\\
				1hnr-	&1&	1&1ecsA	&3 &5 & 1b44D	&1&	1&1ai9A	&1 &1\\
				1bpr-	&1&	1&1dujA	&1&	1&1g31A	&1 &1 & 1auz-	&1&	1\\
				1dax-	&3&	1&1aalA	&1&	1&1i8nA	&1 &21 & 1exg-	&2&	1\\
				1b4rA	&1&	1&1aiw-	&1 &1 & 1hp8-	&1&	1&1qgeE	&1 &1\\
				1a2b-	&1&	1&1xbd-	&7&	1&1kbhA	&1 &1 & 1dc7A	&1&	1\\
				1a2kA	&1&	4&1ab1-	&1&	1&1cqkA	&1 &1 & 1j5eP	&1&	1\\
				1o7bT	&1&	1&1be9A	&1 &1 & 1b4uA	&1&	1&1bdyA	&1 &1\\
				1j7dB	&1&	1&1occE	&1&	1&1cqqA	&1 &1 &1ly7A	&1&	1\\
				1jacA	&1&	1&1hs7A	&1&	3&1ibxB	&1 &1 & 1jajA	&1&	1\\
				1occJ	&1&	2&1qj8A	&1 &1 & 1occH	&1&	6&2sob-	&1 &1\\
				1tafB	&1&	1&1b6q-	&1&	1&1akjD	&2 &1 & 1a4aA	&1&	1\\
				1ha8A	&1&	1&1olgA	&21&	42&1icfA	&1 &1 & 1a4yB	&1&	1\\
				1yuf-	&6&	13&1qjzA	&1 &1 & 1jy2N	&2&	36&1csbA	&2 &1\\
				1dhn-	&1&	1&1cdzA	&1&	12&1qkfA	&1 &1 & 1hbiA	&1&	1\\
				1f7lA	&1&	1&1pcfA	&1&	1&1g84A	&2 &1 & 1jh3A	&1&	1\\
				1cfaA	&1&	1&1kilE	&1 &1 & 1fpzA	&1&	1&1ehxA	&1 &1\\
				1b01A	&3&	13&1hykA	&2&	1&1am9A	&1 &1 & 1a6l-	&1&	1\\
				1dk7A	&1&	1&1hyp-	&1&	4&1amx-	&1 &1 & 1jhcA	&1&	1\\
				1perL &1& 4 & 1qmtA & 1 & 1 & 1cg5B & 1 & 1 & 3lriA & 4 & 1\\
				1fw9A & 17 & 1 & 1gd7A & 1 & 1
				\\
				\hline 
			\end{tabular}
			\begin{tablenotes}
				\item For \hspace{1.5cm} 150 \hspace{1.5cm} 150
				\item Ave  \hspace{1.5cm} 1.67  \hspace{1.5cm} 2.83
				\item First  \hspace{1.5cm} 130 \hspace{1.5cm} 124
			\end{tablenotes}
		\end{threeparttable}
	}
\end{table}

Other force fields such as HR \cite{Rajgaria}, LKF \cite{LKF}, TE13 \cite{TE13}, and HL \cite{HL} have been tested on this set of high-resolution decoys. All these force fields were fundamentally different from each other in their methods of energy estimation. The HR force field is a novel C$_{\alpha}$-C$_{\alpha}$ distance-dependent potential where the interaction distance range 3-9\AA~ is divided into eight bins. Also, this method has the nearest results to our model with a detection limit of 113 native structures. The LKF force field is a C$_{\alpha}$-C$_{\alpha}$ distance-dependent potential where the interaction distance range 3-9\AA~ is divided into eight bins. However, the LKF was not successful in the detection of native structures. Moreover, the TE13 force field is also a distance-dependent (13 bin) force field, but the interaction distance is measured between the geometric centers of the side chain of two interacting residues. The HL force field is contact-based, such that two conditions are needed to consider a pairwise amino acid contact. First, they have to be at least five residues apart from each other, and second, the distance between their non-hydrogen atoms must be less than 4.5 angstroms. The comparison of the energy rankings obtained using different force fields is presented in Table \ref{table3}. As one can see in table \ref{table3}, according to the ranks, LPKPs are the best in identifying the native structures. Assigning rank 1 to the native structure means that the force field is expert at finding the native structure from an array of its nonnative, called decoy structures, configurations. therefore, the LPKPs could increase the percentage of native fold recognition from 75.33 to 86.67.
\begin{table}[h]
	\centering
	\begin{threeparttable}
		\caption{Testing force fields on 150 proteins of the Titan-HRD decoy set.}
		\label{table3}
		\begin{tabular}{p{3cm}ccc}
			\textbf{Method Name} & \textbf{Ave Rank} & \textbf{No of Firsts} & \textbf{Ave RMSD}\\\hline\hline
			LPKP$^{1}$ & 1.67 & 130 (86.67\%) & 2.291\\
			LPKP$^{2}$ & 2.83 & 124 (82.67\%) & 2.291\\
			HR$^{a}$ & 1.87 & 113 (75.33\%) & 0.451\\
			LKF$^{b}$ & 39.45 & 17 (11.33\%) & 1.721\\
			TE13$^{c}$ & 19.94 & 92 (62.16\%) & 0.813\\
			HL$^{d}$ & 44.93 & 70 (46.67\%) & 1.092\\
			\hline
		\end{tabular}
		\begin{tablenotes}
			\item LPKP$^{1}$ and LPKP$^{2}$ are the results of schemes \eqref{5} and \eqref{8}, respectively.
			\item $^{a}$ is extracted from \cite{Rajgaria}.
			\item $^{b}$ is extracted from \cite{LKF}.
			\item $^{c}$ is extracted from \cite{TE13} and TE13 force field was only tested on 148 cases.
			\item $^{d}$ is extracted from \cite{HL}.
		\end{tablenotes}
	\end{threeparttable}
\end{table}

The results of optimization \ref{eq:5} and \ref{eq:8} were tested to calculate the correlation coefficient of each of the proteins in the test Titan-HRD set. We used the below formulation for computing the correlation coefficient:
\begin{equation}
	Corr(d,E) = \frac{1}{N-1} \sum^{N}_{i=1} \frac{S_{d}(X_{i}) - \langle S_{d} \rangle}{\sigma(S_{d})} \frac{E(X_{i}) - \langle E \rangle}{\sigma(E)}.
\end{equation}

\begin{figure}[h]
	\centering
	\includegraphics[width=17cm]{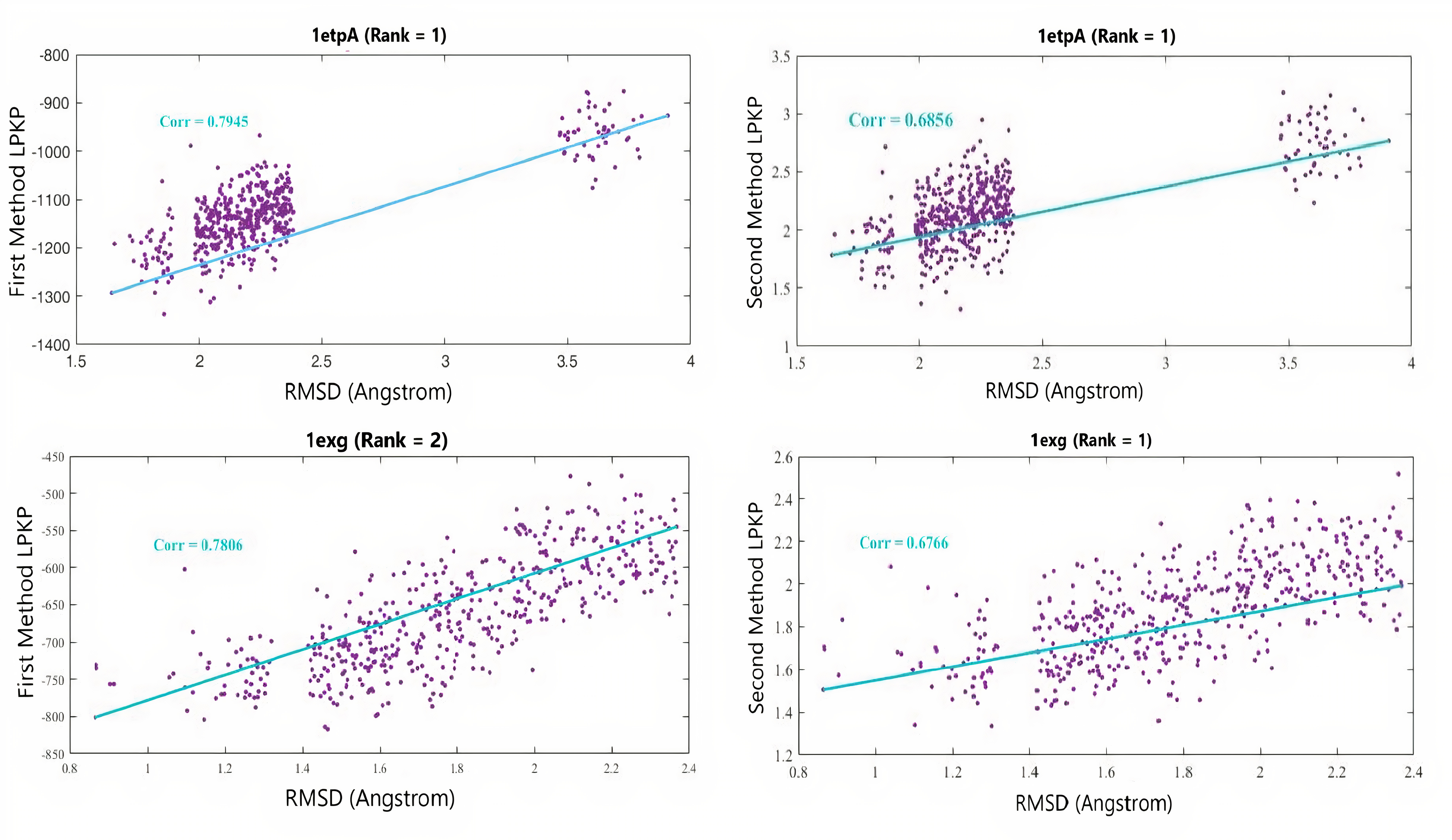}
	\caption{Energy-RMSD plot for two test cases with two methods LPKP$^{1}$ and LPKP$^{2}$.}
	\label{Fig2}
\end{figure}

The $Corr(d,E)$ defined between values of energy $E(X_{i})$ and distance $S_{d}(X_{i})$ for all of the $X_{i}$ decoy structures. Furthermore, $\langle .\rangle$ and $\sigma(.)$ were determined by the mean and standard deviation. The reason for defining the correlation coefficient $Corr(d, E)$ is that it can measure the quality of energy function $E$ concerning the distance $S_{d}$. Initially, by examining the correlation of each protein, we found that scheme \ref{eq:5} had better than \ref{eq:8} in the correlation coefficient LPKP and RMSD. For the presentation of these results, two proteins were randomly chosen in Figure \ref{Fig2} and compared for two optimization methods. The correlation had been approximately 0.7 for the first method and this value reduced to 0.6 in the second.

\begin{figure}[h]
	\centering
	\includegraphics[width=17cm]{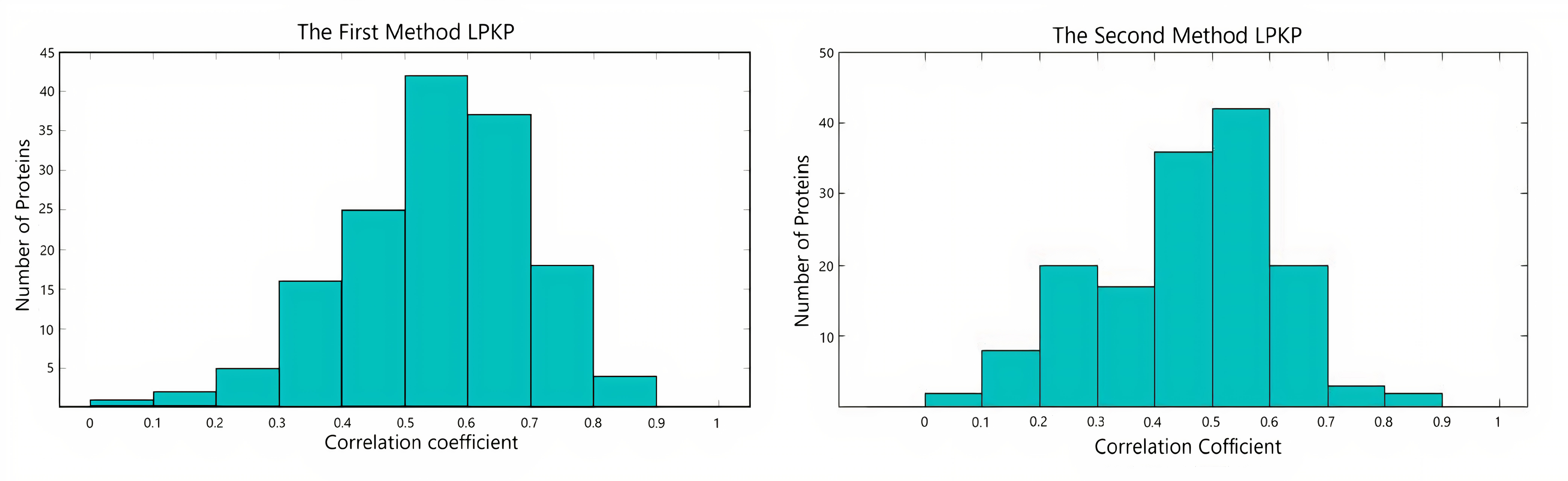}
	\caption{Distribution of correlation coefficient between energy and RMSD for all 150 test cases.}
	\label{Fig3}
\end{figure}

\begin{figure}[h]
	\centering
	\includegraphics[width=10cm]{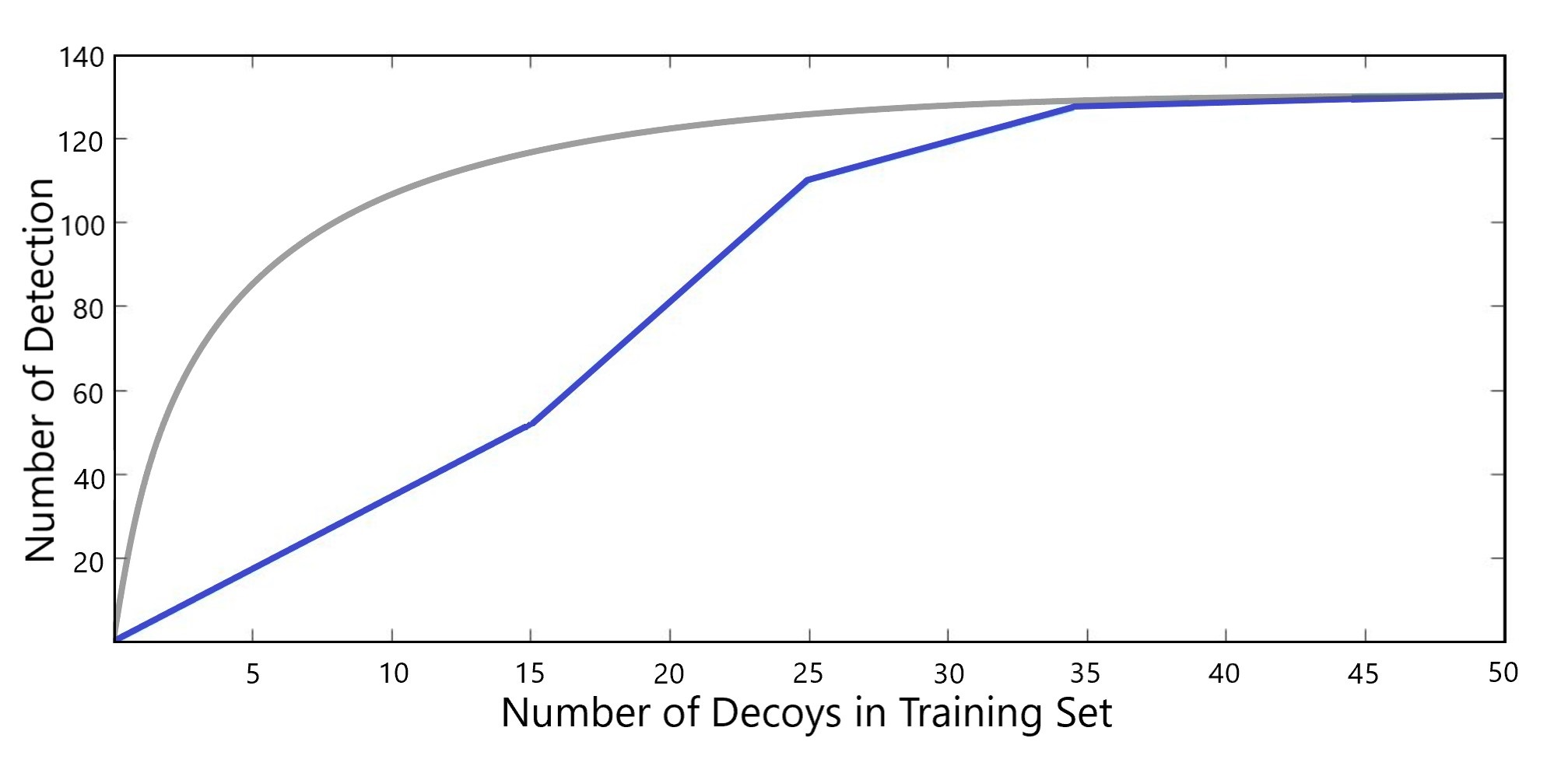}
	\caption{The effect of decoys number in the training set on native structure detection and logarithmic function.}
	\label{Fig4}
\end{figure}

Considering the correlation graph of all proteins in Figure \ref{Fig3}, we found that scheme \ref{eq:5} had a mean correlation of 0.7. Hence, the scheme \ref{eq:5} had more than 80 proteins with a correlation of 0.6. These results indicate a high correlation in LPKP. Furthermore, by comparing the two graphs in Figure \ref{Fig3}, we found that scheme \ref{eq:5} has a higher correlation amount. Also, method \ref{eq:8} had a mean correlation of 0.6 and had more than 70 proteins with a correlation of 0.6.
Originally, the number of decoys and the value of epsilon were considered constant in our training set. Since the LPKP$^{1}$ brought on better results, in the following, the conditions are discussed for this modeling method. Hence, we considered five different values for the number of training set decoys. In this experiment, we included the values of 15, 25, and 35 instead of 45 for the number of decoys in each protein. The results of identifying native structures for each of these values are summarized in Table \ref{table5}. We observed that as the number of decoys increased to 35, our model accuracy increased dramatically. However, this had no significant effect on the rate of detection of the native structure. Furthermore, according to Figure \ref{Fig4}, in each protein, the graph of different values for the number of training set decoys is very similar to the logarithmic function. Considering this diagram, if the number of decoys selected for the training set is more than a certain amount, there is no significant change in the results and only increases the computational volume of the modeling. Therefore, we selected 45 decoy structures per protein for use in the process of developing knowledge-based potentials. In the last step, different values for $\varepsilon$ have been considered in the optimization schemes. We realized the slight changes in the amount of $\varepsilon$ do not create a significant impact on our modeling by examining these values.

\begin{table}[h]
	\centering
	\begin{threeparttable}
		\caption{The effect of the number of decoys in training set on native structure detection.}
		\label{table5}
		\begin{tabular}{cc}
			\textbf{Number of Decoys} & \textbf{No of Firsts}\\\hline\hline
			15 & 52\\
			25 & 110\\
			35 & 128\\
			45 & 130\\
			50 & 130\\
			\hline
		\end{tabular}
		\begin{tablenotes}
			\item Number of decoys for each protein in the training set.
			\item No of firsts is the number of native structures detected in all the test set proteins.
		\end{tablenotes}
	\end{threeparttable}
\end{table}

Although our priority was detecting the native structure among decoys as the first structure, and our proposed scheme was able to present adequate detection regardless of the elimination of many involved parameters and mathematical modeling, for comparison of the LPKPs with more recent methods, the best decoy introduced in \cite{PlosONE} is used, which is demonstrated in Table \ref{table6}.
\begin{table}[h]
	\centering
	\begin{threeparttable}
		\caption{Assessing the best decoys selected by energy functions on Titan-HRD dataset.}
		\label{table6}
		\begin{tabular}{p{2.5cm}ccccc}
			& \textbf{Best} & \textbf{PPD$^{a}$} & \textbf{PPE$^{a}$} & \textbf{LPKP$^{1}$} & \textbf{LPKP$^{2}$}\\
			Titan-HRD & 1.11 & 1.70 & 1.60 & 1.62 & 1.81\\
			\hline
		\end{tabular}
		\begin{tablenotes}
			\item Average value over the test set.
			\item LPKP$^{1}$ and LPKP$^{2}$ are the results of schemes \ref{eq:5} and \ref{eq:8}, respectively.
			\item $^{a}$ is extracted from \cite{PlosONE}.
		\end{tablenotes}
	\end{threeparttable}
\end{table}

\section{Conclusion}
Knowledge-based potentials are developed to derive native structures from their decoy sets. On the other hand, mathematical formulations are designed to increase the accuracy and efficiency of the modeling. In the current study, we have introduced and developed a method based on elementary mathematical methods to display the simplest form of proteins.  the benefits of these methods are easy to implement with high performance and at the same time, with no need for more CPU time compared to popular methods. We constructed two different sets of distance-dependent potential energy functions based on two basic assumptions. At first, we assumed that the energy of each decoy should be more positive than the corresponding native type. The next step assumes that the energy difference and the distance between the two structures are linearly dependent. The RMSD was used to calculate the distance between the decoys and native structures. Each of the potential energy functions has terms of pairwise distances between C$_{\alpha}$-C$_{\alpha}$ and is expressed using the B-spline basis function. We optimized the parameters of the potential function by using two linear programming problems on a large collection of Titan-HRD decoy sets. Furthermore, the obtained results were tested on the remainder of Titan-HRD. We found that the potential functions developed based on Anfinsen’s dogma have more accurate detection than those developed by the root-mean-square deviation of structures. However, both linear programming knowledge-based potentials (LPKP) were successful in recognizing the native structures from an ensemble of high-resolution decoys. In this stage, the high detection was targeted and the LPKP in the first scheme was able to correctly identify 130 native structures out of 150 test cases with an average rank of 1.67. The second LPKP scheme was able to detect 124 native structures with an average rank of 2.83. Although there are methods with very strong predictions in this field, our result indicates that linear programming is a promising method for generating knowledge-based potential functions. All the structures including training and testing Titan-HRD used for this work are available online and can be downloaded from \url{http://titan.princeton.edu/HRDecoys}. Although many influential factors are involved in the production of knowledge-based potentials. Hence in this study, the modeling of these functions with a mathematical approach is discussed. Given that the purpose of modeling is to reduce the parameters involved in the problem as much as possible. Further, our issue was that despite the reduction in computational volume, the designed method was able to score the best structures in those ensembles to produce improved structures. Furthermore, we will attempt to upgrade our approach to an all-atom model in the future.

\bibliographystyle{unsrt}  
\bibliography{main}

\end{document}